\documentclass[MSNbibl,number,citesort,dvips]{arxbj}

\aid{0}
\volume{20}
\issue{4}
\pubyear{2014}
\firstpage{1819}
\lastpage{1844}
\doi{10.3150/13-BEJ543} %kopijuoti is PTS

\makeatletter

\newcommand{\eqref}[1]{(\ref{#1})}
\newcommand{\cb}{\operatorname{CB}}
\newcommand{\cbi}{\operatorname{CBI}}
\newtheorem{thmm}{Theorem}
\newtheorem{lem}[thmm]{Lemma}
\newtheorem{cor}[thmm]{Corollary}
\newtheorem{prop}[thmm]{Proposition}
\newremark{rem}{Remark}[section]
\newremark{remark}{Remark}[section]
\newremark{Example}{Example}[section]
\newproclaim{Def}{Definition}
\makeatother

\begin{document}
\begin{frontmatter}

\title{Local extinction in continuous-state branching processes with
immigration}

\runtitle{Local extinction in CBI processes}

\begin{aug}
%%%% inicialai - be tarpu
\author[a]{\inits{C.}\fnms{Cl\'ement} \snm{Foucart}\thanksref{a}\ead[label=e1]{foucart@math.tu-berlin.de}}
\and
\author[b]{\inits{G.}\fnms{Ger\'onimo} \snm{Uribe Bravo}\corref{}\thanksref{b}\ead[label=e2]{geronimo@matem.unam.mx}}
\address[a]{Institut f\"ur Mathematik, Technische Universit\"at Berlin,
RTG 1845, D-10623, Berlin, Germany.\\ \printead{e1}}
\address[b]{Instituto de Matem\'aticas, Universidad Nacional Aut\'onoma
de M\'exico, \'Area de la Investigaci\'on Cient\'ifica,
Circuito Exterior, Ciudad Universitaria, Coyoac\'an, 04510, M\'exico,
D.F.\\ \printead{e2}}

%%\runauthor{} %% auto
\end{aug}

% HISTORY:
\received{\smonth{12} \syear{2012}}
\revised{\smonth{6} \syear{2013}}

% ABSTRACT
%
\begin{abstract}
%The purpose of this article is to study some processes which have for
%zero-set a random cutout set, as defined by Mandelbrot in 1972. A
%simple but quite rich class of such processes is given by the
%continuous-state branching processes with immigration. We provide
%necessary and sufficient conditions for local extinctions to occur and
%study the properties of their zero sets. This study relies on a useful
%connection between random cutout sets and zero sets through spine
%decomposition. Generalized Ornstein-Uhlenbeck processes valued in $
%known results are recovered. In a second part, we investigate the
%zero-sets of different processes which do not belong into the class of
%continuous branching processes with immigration.
The purpose of this article is to observe that the zero sets of
continuous-state branching processes with immigration (CBI) are
infinitely divisible regenerative sets. Indeed, they can be constructed
by the procedure of random cutouts introduced by Mandelbrot in 1972. We
then show how very precise information about the zero sets of CBI can
be obtained in terms of the branching and immigrating mechanism.
\end{abstract}

% KEYWORDS
% visi is mazosios raides ir pagal abecele
%

\end{frontmatter}

%s1 #&#
\section{Introduction}\label{sec1}
The problem of characterizing the zero set of a real-valued random
process is, in general, not straightforward. When the process is a
one-dimensional diffusion, several methods have been developed to study
their zero sets. When dealing with Markov processes with jumps, the
problem is rather involved and remains highly studied. We refer for
instance to the recent survey of Xiao~\cite{Xiao} where fractal
properties for L\'evy processes and other Markov processes are
discussed. In this paper, we characterize the zero set of
continuous-state branching processes with immigration (here called CBI
processes).

Fundamental results on CBI processes, including their characterization
as the large population limit of Galton--Watson processes with
immigration and the complete determination of the generator, are
obtained by Kawazu and Watanabe \cite{KAW}. Since then, this class of
processes have been studied extensively in several directions.
%For instance, Lambert \cite{Lambert} and Duquesne \cite{Duquesne2}
%have studied (continuum) random trees such whose local times evolve as
%a CBIs, Dawson and Li \cite{DawsonLi06} and \cite{Caballero:2010yg}
%present some constructions of CBI processes in terms of Brownian
%motion and Poisson point processes and Stannat \cite{Stannat} and
%Handa \cite{Handa} study properties of their semigroups and their
%stationary distributions. From the view of applications, we mention
%that the class of CBIs serves also as a generalization of the CIR
%model of mathematical finance for the evolution of interest rates.
%More precisely, the class of CBI processes coincides with the
%so-called affine processes with statespace $[0,\infty)$ (cf.

%The contribution of this article is to reveal some properties of their
%zero sets.

The time evolution of the CBI process in general incorporates two kinds
of dynamics:  reproduction and immigration. Indeed, Kawazu and Watanabe
\cite{KAW} show that the law of a CBI is characterized by the Laplace
exponents $\Psi$ and $\Phi$ of two independent L\'evy processes: a
spectrally positive L\'evy process (which describes the reproduction)
and a subordinator (which describes the immigration). %It is of
%essential importance to take into consideration the effect of
%immigration.

One of consequences of introducing immigration is that zero is not an
absorbing but a reflecting state. We provide necessary and sufficient
conditions for zero to be polar, transient or recurrent. These results
are obtained thanks to a connection between the zero set of a CBI and
the random cutout sets defined by Mandelbrot \cite{Mandelbrot}. Let $Y$
be a CBI process started at zero associated to $(\Psi, \Phi)$. Define
the random set
\[
\mathcal{Z}:=\overline{\{t\geq0; Y_{t}=0\}},
\]
and denote by $v_{s}$ the solution to the differential equation
\[
\frac{\mathrm{d}v_{s}}{\mathrm{d}s}=-\Psi(v_{s}) \quad\mbox{and}\quad v_{0+}=\infty.
\]
Our main result is the following theorem.

%th1 #&#
\begin{thmm}\label{zeroset}
\begin{enumerate}
\item[(i)] $\mathcal{Z}=\{0\}$ if and only if $\int_{0}^{1}\exp
[-\int_{1}^{u}\Phi(v_{s})\,\mathrm{d}s ]\,\mathrm{d}u=\infty$.
\item[(ii)] If $\int_{0}^{1}\exp [-\int_{1}^{u}\Phi
(v_{s})\,\mathrm{d}s
]\,\mathrm{d}u<\infty$, then
\begin{enumerate}[(a)]
\item[(a)] The random set $\mathcal{Z}$ is the closure of the range of a
subordinator with Laplace exponent
\begin{eqnarray*}
L(q) &=& \biggl[\int_{0}^{\infty}\mathrm{e}^{-qt}\exp
\biggl(\int_{t}^{1}\Phi (v_{s})\,\mathrm{d}s
\biggr)\,\mathrm{d}t \biggr]^{-1}
\\
&=& \biggl[\int_{0}^{\infty}\mathrm{e}^{-qt}\exp
\biggl( \int_{v_{1}}^{v_{t}}\frac
{\Phi(u)}{\Psi(u)}\,\mathrm{d}u \biggr)\,\mathrm{d}t
\biggr]^{-1}.
\end{eqnarray*}
\item[(b)] The random set $\mathcal{Z}$ has almost surely a positive
Lebesgue measure if and only if $\int_{\theta}^{\infty}\frac{\Phi
(s)}{\Psi(s)}\,\mathrm{d}s<\infty$ (in that case, we say that the zero set is
\textit{heavy}, otherwise the set is \textit{light}).
%$\int_{0}^{1}\Phi(v_{s})\,\mathrm{d}s<\infty$.
%
\item[(c)] The random set $\mathcal{Z}$ is almost surely the union of
closed nonempty intervals if and only if $\Phi$ is the Laplace
exponent of a compound Poisson process.
\end{enumerate}
\end{enumerate}
\end{thmm}

Under an assumption of regular variation on the ratio $R\dvtx u\mapsto\Phi
(u)/\Psi(u)$, we are able to give more details on the zero set. Loosely
speaking, the index of regularity of the map $R$ at $+\infty$ denoted
by $\rho$ measures the strength of the immigration over the
reproduction. We prove for instance that if $\rho>-1$ then $0$ is
polar. When $\rho=-1$, some new constants are involved in the nature of
the state zero. We also get an upper and a lower bound for the
Hausdorff dimension of the zero set.

The connection made between random covering of the real line and the
zero set of a CBI can be extended to another celebrated class of
processes: the so-called generalized Ornstein--Uhlenbeck processes (OU
processes).
%Indeed, it may be easily observed that the zero set of an OU process
%is also infinitely divisible. However we have not been able to give a
%general cutout construction of their zero sets.
Besides obtaining results to generalized OU processes which are CBI
processes, we provide a characterization of the zero set of
Ornstein--Uhlenbeck processes driven by stable L\'evy processes.
%The paper is organized as follows. In Section \ref{prelim}, we recall
%the definition of the continuous branching processes, some fundamental
%properties and a key result on random covering. In Section
%a corollary of Theorem \ref{zeroset} and a general property about the
%fractal dimensions of the zero set. In Section \ref{analysis}, we
%establish precisely the connection with the random cutout set by
%applying the spine decomposition of a CBI (Theorem \ref{spine}). In
%Section \ref{RV}, we deal with a ratio function $R$ with regular
%variation. Regular variation at $+\infty$ leads us to simple
%conditions on the upper and lower indices at $+\infty$ of $\Phi$ and $
%to conditions on the indices of $\Phi$ and $\Psi$ at $0$ for
%transience or recurrence. Finally in Section \ref{OU}, we study the
%zero set of the Ornstein Uhlenbeck process driven by a stable L\'evy
%process.
%We begin by recalling some fundamental definitions and properties
%about continuous-state branching processes with immigration. We refer
%for instance the reader to Chapter 10 of Kyprianou \cite{MR2250061}
%and Chapter 3 of Li \cite{Li} for a comprehensive survey on
%continuous-state branching processes. Then, we recall some theorems
%about random covering of the real line due to Mandelbrot
%random cutout set is at the core of the present work.

The paper is organized as follows.
Section~\ref{prelim} contains preliminaries on CBI processes and
Mandelbrot's random cutout construction.
Section~\ref{Zeroset} contains the first implications of Theorem \ref
{zeroset} to properties of the zero set of CBI processes as polarity,
transience, recurrence, and Box counting and Hausdorff dimensions.
In Section~\ref{analysis}, we link the zero set of a CBI process to
random cutouts through the spine decomposition and prove Theorem \ref
{zeroset}. We then particularize to the case when the ratio of the
immigration and branching mechanisms of our CBI process is regularly
varying in Section~\ref{RV}. Finally, Section~\ref{OU} is devoted to a
characterization of the zero set of an Ornstein--Uhlenbeck process
driven by a stable L\'evy process.
%s2 #&#
\section{Preliminaries}\label{prelim}
%s2.1 #&#
\subsection{Continuous-state branching processes} \label{cb}
Let $x\in\mathbb{R}_{+}$.
A Markov process $X(x)=(X_t(x),t\geq0)$, where $X_0(x)=x$ is called a
branching process in continuous time and continuous space ($\cb$ for
short) if it satisfies the following property: For any $y\in\mathbb{R}_{+}$
\[
X(x+y)\stackrel{d} {=}X(x)+\tilde X(y),
\]
where $\tilde X(y)$ is an independent copy of $X(y)$.

Let $X(x)$ be a $\cb$ issued from $x$. There exists a unique triplet
$(d, \sigma, \nu)$ with $d, \sigma\geq0$, and $\nu$ a~measure carried
on $\mathbb{R}_{+}$ satisfying
\[
\int_{0}^{\infty} \bigl(1\wedge x^{2}
\bigr) \nu(\mathrm{d}x)<\infty
\]
such that the Laplace transform of the one-dimensional distribution of
$X_t(x)$ is given by
\[
\mathbb{E} \bigl[\mathrm{e}^{-\lambda X_{t}(x)} \bigr]=\exp \bigl(-xv_{t}(\lambda)
\bigr),
\]
where the map $t\mapsto v_{t}(\lambda)$ is the solution to the
differential equation
\[\frac{\partial}{\partial t}v_{t}(\lambda)=-\Psi
\bigl(v_{t}(\lambda)\bigr),\qquad v_{0}(\lambda)=\lambda
\]
with
%
%e1 #&#
%
\begin{equation}
\label{Levykhintchine} \Psi(q)=\frac{\sigma
^{2}}{2}q^{2}+{d}q+\int
_{0}^{\infty} \bigl(\mathrm{e}^{-qx}-1+qx1_{x\leq
1}
\bigr)\nu(\mathrm{d}x).
\end{equation}
The process is said to be critical, subcritical or supercritical
according as $\Psi'(0+)=0$, $\Psi'(0+)>0$ or $\Psi'(0+)<0$.
For any $x\geq0$, define the extinction time of the $\operatorname{CB} (X_{t}(x),
t\geq0)$ by
\[
\zeta:=\inf \bigl\{t\geq0; X_{t}(x)=0 \bigr\}.
\]
Recall the following result regarding the distribution of the
extinction time (see, e.g., Theorems~3.5 and~3.8, pages 59--60 of
\cite{Li}):
\[
\mathbb{P}_{x}[\zeta\leq t]=\exp{(-xv_{t})},
\]
where $v_{t}={ \lim}
_{\lambda\rightarrow\infty}v_{t}(\lambda
)$. Also, recall \textit{Grey's theorem}.
%
%th2 #&#
\begin{thmm}[(Grey \cite{Grey})]\label{Grey} The $\operatorname{CB} (X_{t}, t\geq0)$ is
absorbed in $0$ with positive probability if and only if there exists
$\theta>0$ such that $\Psi(z)>0$ for $z\geq\theta$ and
%
%e2 #&#
%
\begin{equation}
\label{AbsorptionHypothesis} \int_{\theta}^{\infty}
\frac{\mathrm{d}q}{\Psi(q)}<\infty.
\end{equation}
Under that integrability condition (called Grey's condition), the real
number $v:={\lim\downarrow}_{t \rightarrow\infty} v_{t}\in
[0, \infty[$ is the largest root of the equation $\Psi(x)=0$ and
\[
\mathbb{P}_{x}[\zeta<\infty]=\exp{(-xv)}.
\]
In the (sub)critical case, $v=0$ and the $\operatorname{CB}(\Psi)$ is absorbed at zero
almost surely. In the supercritical case, $v>0$.
\end{thmm}
%
%Under Grey's condition, the map $(v_{s}, s\geq0)$ is the minimal
%solution to the differential equation (see Corollary 3.11 of
%v_{0+}=\infty. \end{equation}
%s2.2 #&#
\subsection{Continuous-state branching processes with immigration}
\label{CBI}
A continuous-state branching processes with immigration (CBI for short)
started at $x$ is a Markov process $Y(x)=(Y_{t}(x),t\geq0)$ satisfying
the following property: for any $y\in\mathbb{R}_{+}$
\[
Y(x+y)\stackrel{d} {=}Y(x)+ X(y),
\]
where $X(y)$ is an independent $\cb$ with mechanism $\Psi$.
Any CBI process is characterized by two functions of the variable
$q\geq0$:
\begin{eqnarray*}
\Psi(q)&=&{d}q+\frac{1}{2}\sigma^{2}q^{2}+\int
_{0}^{\infty
} \bigl(\mathrm{e}^{-qu}-1+qu1_{u\in(0,1)}
\bigr)\nu_{1}(\mathrm{d}u),
\\
\Phi(q)&=&\beta q+\int_{0}^{\infty}
\bigl(1-\mathrm{e}^{-qu} \bigr)\nu_{0}(\mathrm{d}u),
\end{eqnarray*}
where $\sigma^{2}, \beta\geq0$ and $\nu_{0}$, $\nu_{1}$ are two L\'
evy measures such that
\[
\int_{0}^{\infty} (1\wedge u) \nu_{0}(\mathrm{d}u)<
\infty\quad\mbox{and}\quad \int_{0}^{\infty} \bigl(1\wedge
u^{2} \bigr) \nu_{1}(\mathrm{d}u)<\infty.
\]
The measure $\nu_{1}$ is the L\'evy measure of a spectrally positive
L\'
evy process which characterizes the reproduction. The measure $\nu_{0}$
characterizes the jumps of the subordinator that describes the arrival
of immigrants in the population. The nonnegative constants $\sigma^2$
and $\beta$ correspond, respectively, to the continuous reproduction and
the continuous immigration. Let $\mathbb{P}_{x}$ be the law of a
$\operatorname{CBI}(Y_{t}(x), t\geq0)$ started at $x$, and denote by $\mathbb{E}_{x}$
the associated expectation. The law of $(Y_{t}(x),t\geq0)$ can then be
characterized by the Laplace transform of its marginal as follows: for
every $q>0$ and $x\in\mathbb{R}_{+}$,
\[
\mathbb{E}_{x} \bigl[\mathrm{e}^{-qY_{t}} \bigr]=\exp
\biggl(-xv_{t}(q)-\int_{0}^{t}\Phi
\bigl(v_{s}(q) \bigr)\,\mathrm{d}s \biggr),
\]
where
\[
v_{t}(q)=q-\int_{0}^{t}\Psi
\bigl(v_{s}(q) \bigr)\,\mathrm{d}s.
\]

The pair $(\Psi, \Phi)$ is known as the branching and immigration
mechanisms. A CBI process $(Y_{t}, t\geq0)$ is said to be conservative
if for every $t>0$ and $x\in[0,\infty[, \mathbb{P}_{x}[Y_{t}<\infty
]=1$. A result of Kawazu and Watanabe \cite{KAW} states that $(Y_{t},
t\geq0)$ is conservative if and only if for every $\varepsilon>0$,
\[
\int_{0}^{\varepsilon}\frac{1}{|\Psi(q)|}\,\mathrm{d}q=\infty.
\]
See Theorem 10.3 in Kyprianou \cite{MR2250061} for a proof.

Contrary to the simpler setting of continuous-state branching processes
without immigration, nondegenerate stationary laws may appear. The
following theorem provides a necessary and sufficient condition for the
CBI to have a stationary distribution. Properties of the support of its
stationary law have been studied by Keller-Ressel and Mijatovi\'c in
\cite{Keller}.
%
%th3 #&#
\begin{thmm}[(Pinsky \cite{Pinsky}, Li \cite{Li})]\label{stationary}
The $\operatorname{CBI} (\Psi, \Phi)$ process has a stationary law if and only if
\[
\Psi'(0+)\geq0\quad \mbox{and}\quad\int_{0}^{\theta}
\frac{\Phi(u)}{\Psi
(u)}\,\mathrm{d}u<\infty.
\]
If the process is subcritical ($\Psi'(0+)>0$), the convergence of this
integral is equivalent to the following log-condition
\[
\int_{x\geq1}\log(x)\nu_{1}(\mathrm{d}x)<\infty.
\]
\end{thmm}

The main ingredient of the proofs provided in this work relies on a
connection between the zero set of the CBI and a particular random set,
called random cutout set.
%Many processes have zero set corresponding with random cutout. We
%shall investigate in this article several examples. The easiest will
%be the CBI processes.
%s2.3 #&#
\subsection{Random covering of the real half-line}\label{cutout}
We recall here the definition of a random cutout set studied first by
Mandelbrot \cite{Mandelbrot}.
%The notion of random cutout set is crucial in our study.
For a textbook presentation of the main theorems regarding random
cutouts, we refer to the course of Bertoin \cite{subordinatorbertoin}.
The main results we shall use in this paper may be found in. Consider a
$\sigma$-finite measure $\mu$ on $\mathbb{R}_{+}$ which is finite on
compact subsets of $(0,\infty)$. Denote its tail $\mu([x,\infty[)$ by
$\bar{\mu}(x)$ for all $x\in\mathbb{R}_{+}$. Let $\mathcal{N}$ be a
Poisson point process on $\mathbb{R}_{+}\times\mathbb{R}_{+}$ with
intensity $\mathrm{d}t\otimes\mu$. Denote by $(t_{i}, x_{i})_{i\in\mathcal
{I}}$ its atoms:
\[
\mathcal{N}=\sum_{i\in\mathcal{I}}\delta_{(t_{i}, x_{i})}.
\]
We recover the half line by the intervals $]t_{i},t_{i}+x_{i}[$, $i\in
\mathcal{I}$. The set of uncovered point is
\[
\mathcal{R}= [0, \infty)- \bigcup_{i\in I}
\,]t_{i},t_{i}+x_{i}[.
\]
The measure $\mu$ is called the cutting measure. A random set
$\mathcal
{R}$ is said to be a random cutout set if it is obtained by such a
construction. Any random cutout set is a regenerative set (i.e., the
closure of the range of a subordinator).
%
%th4 #&#
\begin{thmm}[(Fitzsimmons, Fristedt and Shepp \cite{MR799145})] \label
{recovering}
\[
\mbox{If }\int_{0}^{1}\exp \biggl( \int
_{t}^{1}\bar{\mu }(s)\,\mathrm{d}s \biggr)\,\mathrm{d}t=\infty\mbox{ then
} \mathcal{R}=\{0\} \mbox{ a.s.}
\]
Otherwise $\mathcal{R}$ is the closure of the image of a subordinator
with Laplace exponent $L$ given by
\[
\frac{1}{L(q)}=\int_{0}^{\infty}\mathrm{e}^{-qt}
\exp \biggl(\int_{t}^{1}\bar {\mu }(s)\,\mathrm{d}s \biggr)\,\mathrm{d}t.
\]
Moreover if $\int_{1}^{\infty} \exp (\int_{t}^{1}\bar{\mu
}(s)\,\mathrm{d}s )\,\mathrm{d}t=\infty$, then the set $\mathcal{R}$ is unbounded.
\end{thmm}
These random cutout sets have been intensively studied. As mentioned in
the \hyperref[sec1]{Introduction}, we shall prove that the zero set of a CBI process is
a random cutout set. Along the article, we shall use results on their
geometry and refer the reader to the monography of Bertoin \cite
{subordinatorbertoin}. For instance, we mention that Theorem \ref
{recovering} matches with Theorem 7.2 of \cite{subordinatorbertoin}.

We end this section with the notion of infinitely divisibility for a
regenerative set.
%As we shall see, this notion is underlying in all our work.
A regenerative set $\mathcal{R}$ is said to be infinitely divisible if
for any $n\geq1$, there exist $n$ independent identically distributed
regenerative sets $(\mathcal{R}^{i}, 1\leq i \leq n)$ such that
\[
\mathcal{R}\stackrel{\mathrm{law}} {=}\bigcap_{i=1}^{n}
\mathcal{R}^{i}.
\]
Any random-cutout set is infinitely divisible (cf. \cite{MR799145}, Theorem
3). To the best of our knowledge, the converse is still an
open question (cf. \cite{MR2132405}, Open problem 2.24, page 334).
%s3 #&#
\section{The zero set of continuous-state branching processes with
immigration} \label{Zeroset}
%As recalled in Subsection \ref{cb},
In this section, we give some basic consequences of Theorem \ref
{zeroset}. The latter theorem is proved in Section~\ref{analysis}. We
shall focus on a $\operatorname{CBI}(\Psi, \Phi)$ started at $0$ such that $\Psi$
satisfies Grey's condition. As we shall see in Section~\ref{analysis},
this is of no relevance for the study of the zero set. By a slight
abuse of notation, we denote by $(Y_{t}, t\geq0)$ the process
$(Y_{t}(0), t\geq0)$. Define the random set
\[
\mathcal{Z}:=\overline{\{t\geq0; Y_{t}=0\}}.
\]

%pr5 #&#
\begin{prop} The random set $\mathcal{Z}$ is infinitely divisible.
\end{prop}
\begin{pf} Let $n\geq1$. Consider $n$ independent CBI processes
$(Y^{i}; 1\leq i\leq n)$ started at $0$ with branching mechanism $\Psi$
and immigrating mechanism $\frac{1}{n}\Phi$. We have clearly the
following equalities in law
\[
(Y_{t}, t\geq0)\stackrel{\mathrm{law}} {=} \Biggl(\sum
_{i=1}^{n}Y^{i}_{t}, t\geq0
\Biggr),
\]
and therefore
\[
\mathcal{Z}\stackrel{\mathrm{law}} {=}\bigcap_{i=1}^{n}
\overline{ \bigl\{t\geq0; Y^{i}_{t}=0 \bigr\}}.%\qedhere
\]
\upqed\end{pf}
%
%Our main result is the following theorem.
%%\item[i)] $\mathcal{Z}=\{0\}$ if and only if $\int_{\theta}^{\infty}
%%\item[ii)] If $\int_{\theta}^{\infty}\exp\left[\int_{\theta}^{z}\frac{
%a subordinator with Laplace exponent
%%\begin{center}
%L(q)
%&=\left[\int_{0}^{\infty}e^{-qt}\exp\left(\int_{t}^{1}\Phi(v_{s})\,\mathrm{d}s
%%\end{center}
%Lebesgue measure if and only if $\int_{\theta}^{\infty}\frac{\Phi(s)}{
%%$\int_{0}^{1}\Phi(v_{s})\,\mathrm{d}s<\infty$.
%closed non-empty intervals if and only if $\Phi$ is the Laplace
%exponent of a compound Poisson process.
%In order to study the transience/recurrence of the state zero for
%$(Y_{t}, t \geq0)$, we appeal to the renewal measure (up to a
%multiplicative factor) of the underlying subordinator. It is
%absolutely continuous and its density is given by $$u(t)=\exp\left(
If the set $\mathcal{Z}$ is not reduced to $\{0\}$ and bounded
(resp., unbounded), the state $0$ is transient (resp.,
recurrent). Recall that the state $0$ is said to be polar if
\[
\mathbb{P}_{x}[\exists t\geq0, Y_{t}=0]=0
\]
for any $x\neq0$.
%
%re3.1 #&#
\begin{rem}
If the process is supercritical, then clearly $0$ is polar or
transient. This can be easily observed using Theorem \ref{zeroset}.
\end{rem}
%
%the underlying subordinator is killed. This is equivalent to $L(0)>0$,
%and therefore to
%$$\int_{1}^{\infty}\exp\left(\int_{t}^{1}\Phi(v_{s})\,\mathrm{d}s\right)\,\mathrm{d}t<
%Since the process is supercritical $(v_{t}, t>0)$ converges as $t\to
%map $\Phi$ is continuous, then $(\Phi(v_{t}), t\geq1)$ also converges
%and is bounded. Denote by $C$ a lower bound of $\Phi(v_{t})$. We have $
The next result is a corollary of Theorem \ref{zeroset}. Since the
supercritical case is plain, we focus on (sub)critical reproduction
mechanism when studying the recurrence and transience of $0$
(statements (ii) and (iii) below).
%
%co6 #&#
\begin{cor} \label{recurrence} Let $\theta>0$ such that $\Psi(u)>0$ for
all $u\geq\theta$.
\begin{longlist}[(iii)]
\item[(i)] The state $0$ is polar if and only if
\[
\int_{\theta}^{\infty}\exp \biggl[\int_{\theta}^{z}
\frac
{\Phi(u)}{\Psi(u)}\,\mathrm{d}u \biggr]\frac{1}{\Psi(z)}\,\mathrm{d}z=\infty.
\]
Assume further that $\Psi$ is (sub)critical, the state $0$ is
\item[(ii)] transient if and only if
\[
\int_{\theta}^{\infty}\exp \biggl[\int_{\theta}^{z}
\frac{\Phi
(u)}{\Psi
(u)}\,\mathrm{d}u \biggr]\frac{1}{\Psi(z)}\,\mathrm{d}z<\infty\quad\mbox{and}\quad \int
_{0}^{\theta}\exp \biggl[ -\int_{x}^{\theta}
\frac{\Phi(s)}{\Psi(s)}\,\mathrm{d}s \biggr] \frac{\mathrm{d}x}{\Psi(x)}<\infty.
\]
\item[(iii)] recurrent if and only if
\[
\int_{\theta}^{\infty}\exp \biggl[\int_{\theta}^{z}
\frac{\Phi
(u)}{\Psi
(u)}\,\mathrm{d}u \biggr]\frac{1}{\Psi(z)}\,\mathrm{d}z<\infty\quad\mbox{and}\quad \int
_{0}^{\theta}\exp \biggl[ -\int_{x}^{\theta}
\frac{\Phi
(s)}{\Psi(s)} \,\mathrm{d}s \biggr] \frac{\mathrm{d}x}{\Psi(x)}=\infty.
\]
\end{longlist}
\end{cor}
%
%re3.2 #&#
\begin{rem}
We may compare the integral conditions with those for the
continuous-state branching process without immigration. The first
statement has to be compared with Grey's condition.
%Note that it differs from Corollary 1.1 in Ogura \cite{MR0279902}
%which gives a necessary and sufficient condition for a fixed point to
%be a zero of a CBI process.
The probability of blow-up of the CB (and CBI) depends on the behaviour
of $\Psi$ at~$0$. Assume that the process is conservative, which holds
if and only if $\int_{0}^{\theta}\frac{\mathrm{d}x}{\Psi(x)}=\infty$. If
moreover $\int_{0}^{\theta}\frac{\Phi(s)}{\Psi(s)}\,\mathrm{d}s<\infty$,
then the
second condition in statement (iii) is always satisfied. Note that if
the process is not conservative, obviously the state $0$ is either
polar or transient.
\end{rem}
\begin{pf*}{Proof of Corollary \ref{recurrence}}
Statements of the
corollary are easily obtained by substitution in the integrals
appearing in Theorem \ref{zeroset}. Recall that $v_{0+}=\infty$, the
first statement is equivalent to the statement (i) of Theorem \ref
{zeroset} by considering $u=v_{s}$. The transience or the recurrence of
zero hold, respectively, if the set $\mathcal{Z}$ is bounded or
unbounded. As already explained, this is equivalent for the underlying
subordinator to be killed or not. If the state $0$ is not polar, then
the zero set is bounded if and only if $L(0)>0$. By the same
substitution in the integral condition of Theorem~\ref{zeroset}, we get
$\int_{1}^{\infty}\Phi(v_{t})\,\mathrm{d}t=\int_{v}^{v_{1}}\frac{\Phi
(u)}{\Psi
(u)}\,\mathrm{d}u$. The constant $v$ is the limit of $(v_{t}, t\geq0)$ when $t$
goes to $\infty$ and equals to $0$ since we focus on the (sub)critical
case. The constant $v_{1}$ does not play any role here and may be
replaced by $\theta$.
\end{pf*}
Recall the definition of the lower and upper box-counting dimension.
See, for instance, Section~3 of \cite{Xiao} or Chapter~5 of \cite
{subordinatorbertoin}. For every nonempty bounded subset $E$ of
$\mathbb{R}_{+}$, let $N_{\varepsilon}(E)$ be the smallest number of
intervals of length $\varepsilon$ needed to cover $E$. The upper and lower
box-counting dimension of $E$ are defined as
\[
\overline{\operatorname{dim}}(E):=\mathop{\limsup }_{\varepsilon\rightarrow0}\frac{\log(N_{\varepsilon}(E))}{\log
(1/\varepsilon)}
\]
and
\[
\underline{\operatorname{dim}}(E):=\mathop{\liminf }_{\varepsilon\rightarrow0}\frac{\log(N_{\varepsilon}(E))}{\log
(1/\varepsilon)}.
\]
We provide now a last general result on the zero set. The proof is
postponed at the end of Section~\ref{analysis}.
%
%le7 #&#
\begin{lem}\label{Hausdorff}
The random set $\mathcal{Z}$ has the following upper and lower
box-counting dimensions. $\mbox{For every } t>0$,
\[
\overline{\operatorname{dim}} \bigl(\mathcal{Z}\cap[0,t] \bigr)=1-\mathop{\liminf
}_{u\rightarrow
0} \ \frac{1}{\log({1}/{u})}\int_{v_{1}}^{v_{u}}
\frac{\Phi
(s)}{\Psi(s)}\,\mathrm{d}s \qquad\mbox{a.s.}
\]
and
\[
\underline{\operatorname{dim}} \bigl(\mathcal{Z}\cap[0,t] \bigr)=1-\mathop{\limsup
}_{u\rightarrow
0} \ \frac{1}{\log({1}/{u})}\int_{v_{1}}^{v_{u}}
\frac{\Phi
(s)}{\Psi(s)}\,\mathrm{d}s\qquad \mbox{ a.s.}
\]
Moreover if $\mathcal{Z}$ is bounded almost surely, then the law of the
last zero of $(Y_{t}, t\geq0)$,
\[
g_{\infty}:=\sup\{s\geq0; s\in\mathcal{Z}\}
\]
is given by
\[
\mathbb{P}[g_{\infty}\in \mathrm{d}t]=k^{-1}\exp \biggl(\int
_{t}^{1}\Phi (v_{s})\,\mathrm{d}s
\biggr)\,\mathrm{d}t=k^{-1}\exp \biggl(\int_{v_{1}}^{v_{t}}
\frac{\Phi
(u)}{\Psi(u)}\,\mathrm{d}u \biggr)\,\mathrm{d}t,
\]
with $k$ the renormalization constant.
\end{lem}
%
%re3.3 #&#
\begin{rem} The lower box-counting dimension %(denoted by $\underline{
and the Hausdorff dimension (denoted by $\operatorname{dim}_{H}$) of a
regenerative set coincides almost surely (see Corollary 5.3 of \cite
{subordinatorbertoin}).
\end{rem}

%s4 #&#
\section{Analysis of the zero set: Spine decomposition and random
covering}\label{analysis}

The objective of this section is to recall a construction of CBI
processes which will then be used to prove Theorem \ref{zeroset}.

A $\operatorname{CB}(\Psi)$ process reaches zero with positive probability if and only
if the branching mechanism $\Psi$ satisfies Grey's condition.
Intuitively, a $\operatorname{CBI}(\Psi, \Phi)$ cannot touch zero, unless the
corresponding $\operatorname{CB}(\Psi)$ can reach zero. We begin by establishing
rigorously this idea. Let $x>0$; the branching property stated in
Section~\ref{CBI} provides that the $\operatorname{CBI}(\Psi, \Phi)$ process $Y(x)$
satisfies
\[
Y(x)\stackrel{\mathrm{law}} {=}Y(0)+X(x)
\]
for a $\operatorname{CB}(\Psi)$ $X(x)$ that starts at $x$. If $x>0$ and $Y(x)$ reaches
zero with positive probability, then plainly the $\operatorname{CB}(\Psi)$ $X(x)$
reaches also zero with positive probability. Theorem \ref{Grey} ensures
then that $\Psi$ verifies Grey's condition. In order to handle the case
$x=0$, we note that $\Phi\neq0$ implies that
%$\imf{v_t}{\lambda}>0$
%
\[
\mathbb{E} \bigl( \mathrm{e}^{-\lambda Y_t(0)} \bigr)=\mathrm{e}^{-\int_0^t \Phi
(v_s(\lambda)) \,\mathrm{d}s}<1
\]
and so $Y_t(0)>0$ with positive probability for any $t>0$. If $Y(0)$
reaches zero with positive probability, then it does so after time $t$
for some $t>0$. Applying the Markov property at time~$t$ on the set
$Y_t(0)>0$, we are reduced to the previous case and conclude that
Grey's condition holds.
Thus, we see that imposing Grey's condition on the branching mechanism
merely rules out a trivial case in which the zero-set of $Y(0)$ is
almost surely~$\{0\}$.

%Notice that we also get that with positive probability $\tau-\sigma$
%is a zero of the process $Y'$ (which has the same law as $Y$) and $
%other words, roughly speaking, the process hits zero for arbitrary
%small times. We shall provide more precise arguments in the sequel.
%%Notice that we also get that $\tau-t$ is a zero of the process $Y'$
%(which has the same law as $Y$) with positive probability. In other
%words, roughly speaking, the zero $\tau$ is not isolated in law. We
%shall provide more precise statements in the sequel.
In order to prove Theorem \ref{zeroset}, we now recall a Poissonian
decomposition of CBI processes, also called the spine decomposition,
which is obvious in the simpler setting of Galton--Watson processes in
immigration and indeed has been considered by a number of authors in
the continuous setting. For example, \cite{MR656509} use it to obtain a
decomposition of Bessel bridges or in \cite{MR1249698} it allows
representations of superprocesses conditioned on nonextinction, a work
which has continued in~\cite{MR2006204,MR2108152}.

The spine decomposition of CBI processes is found in Section~2.1 of
\cite{spinedecomposition}
or as a particular case of the construction of Section~2.2 of \cite
{Abraham}. In the case of CB processes with stable reproduction
mechanism conditioned on nonextinction, which are self-similar CBI
processes, the spine decomposition is found in \cite{MR2582640}. The
spine decomposition is based on the $\mathbb{N}$-measures constructed
for superprocesses in \cite{MR2092876} and specialized to the case of
CBI processes in Theorem 1.1 of \cite{spinedecomposition}. We now give
a streamlined exposition of the construction of this specialized
$\mathbb{N}$-measure, based on \cite{MR656509}, by assuming Grey's condition.

Starting in \cite{MR656509}, a Poisson process representation for
(continuous) CBI processes has been achieved by means of a $\sigma
$-finite measure which can be understood as the excursion law of a $\cb
(\Psi)$ although, in general, one is not able to concatenate excursions
to obtain a recurrent extension of a $\cb(\Psi)$ (cf. \cite{MR656509},
page~440).
%See for instance Subsection 3.1 of \cite{Lambert}.

For a $\cb ( \Psi) $, we have
\[
\mathbb{E}_x^{\Psi} \bigl( \mathrm{e}^{-\lambda X_t} \bigr)
=\mathrm{e}^{-x v_t ( \lambda) }
\]
for any $x\geq0$. Hence, $\lambda\mapsto v_t ( \lambda) $ is the
Laplace exponent of a subordinator, under \eqref{AbsorptionHypothesis},
it is a driftless subordinator since
\[
\mathop{\lim}_{\lambda\to\infty}\frac{v_t ( \lambda) }{\lambda}=0.
\]
See, for instance, Corollary 3.11 page 61 of \cite{Li}. (This is the
fundamental simplification in comparison with \cite{spinedecomposition}.)
Let $( P_t,t\geq0) $ be the semigroup of the $\cb ( \Psi) $ and
$\eta_t$ stand for the L\'evy measure of the Laplace exponent $v_t$ so that
\[
v_t ( \lambda) =\int_0^\infty
\bigl( 1-\mathrm{e}^{-\lambda x} \bigr) \eta_t (\mathrm{d}x) .
\]
By the composition property $v_t\circ v_s(\lambda)=v_{t+s}(\lambda)$
and so
\begin{eqnarray*}
\lambda\int \mathrm{e}^{-\lambda y}\eta_{t+s}  (y,\infty)  \,\mathrm{d}y
&=&
\int1-\mathrm{e}^{-\lambda x} \eta_{t+s} (\mathrm{d}x)
=v_{t+s} ( \lambda) %\\&=\imf{v_t}{\imf{v_s}{\lambda}}
\\
&=&\int1-\mathrm{e}^{-xv_s ( \lambda) } \eta_t (\mathrm{d}x) %\\&=\int\imf{\mathbb{E}_x}{1-\mathrm{e}^{-\lambda X_s}} \imf{\eta_t}{dx}
\\
&=&\lambda\int \mathrm{e}^{-\lambda y }\int\mathbb{P}_x (
X_s>y) \eta_t (\mathrm{d}x) \,\mathrm{d}y.
\end{eqnarray*}
We then see that $\eta_{t+s}$ and $\eta_tP_s$ have the same tails and
that therefore, $( \eta_t,t>0) $ is an entrance law for the
semigroup $P_t$. Hence, we may consider the $\sigma$-finite measure $Q$
on the space of c\`adl\`ag excursions starting and ending at zero
characterized by its finite-dimensional marginals: for $0<t_1<\cdots<t_n$
\[
Q ( X_{t_1}\in \mathrm{d}x_1,\ldots,X_{t_n}\in
\mathrm{d}x_n) =\mathbf{1}_{x_1,\ldots
,x_n>0}\eta_{t_1} (\mathrm{d}x_1) P_{t_2-t_1}^\Psi ( x_1,\mathrm{d}x_2)
\cdots P_{t_n-t_{n-1}}^\Psi ( x_{n-1},\mathrm{d}x_n)
.
\]
This can be thought of as an excursion law for the $\cb ( \Psi) $.
It is not a trivial point that $X_t\to0$ as $t\to0+$ $Q$-almost
everywhere. Pitman and Yor argue, in the case of diffusions, by stating
a William's type decomposition of the measure $Q$: on the set where $X$
reaches a height $>x$, the measure $Q$ is proportional to the
probability measure which concatenates the law of a $\cb(\Psi)$
conditioned to stay positive (started at zero) until the process
reaches a height $>x$ with a trajectory of a $\cb(\Psi)$ until it
reaches zero. Although not explicitly stated, this point of view is the
basis for the proof given in \cite{2012arXiv1202.3223L}, Section~2.4,
that $X_{0+}=0$ under $Q$. Also, under $Q$, $X$ is Markovian and with
the semigroup of the $\cb ( \Psi) $. To characterize the image of
the length of the excursions under $Q$, say $\zeta$, notice that
\begin{eqnarray*}
Q ( \zeta>t) &=&\mathop{\lim}_{h\to0+}Q ( \zeta> t+h)
=\mathop{\lim}_{h\to0+}\int\mathbb{P}_x^\Psi (
\zeta> t) \eta_h (\mathrm{d}x)
\\
&=&\mathop{\lim}_{h\to0+}\int \bigl( 1-\mathrm{e}^{-xv_t} \bigr)
\eta_h (\mathrm{d}x)
=\mathop{\lim}_{h\to0+}v_h ( v_t)
=v_t.
\end{eqnarray*}
Now use $\Phi$ and $Q$ to construct the $\sigma$-finite measure
$\mathbb
{N}$ on excursion space by means of
%$$
%$$which will be the intensity of the Poisson point process
%$$
%
\[
\mathbb{N}=\beta Q+\int\nu_0 (\mathrm{d}x)
\mathbb{P}_x^\Psi,
\]
which will be the intensity of the Poisson point process
\[
\Theta=\sum_{i\in\mathcal{I}} \delta_{( t_i,X^{i}) }.
\]
Finally, let
\[
Y_t=\sum_{t_i\leq t}X^{i}_{t-t_i}.
\]

%th8 #&#
\begin{thmm}\label{spine}
$Y$ is a $\cbi ( \Psi,\Phi) $ which starts at $0$.
\end{thmm}
This is what we refer to as a spine decomposition of a CBI process.
\begin{pf*}{Proof of Theorem \ref{spine}}
We start with a preliminary computation which uses the exponential
formula for Poisson point processes:
\begin{eqnarray*}
\mathbb{E} \bigl(\mathrm{e}^{-\lambda Y_t} \bigr) &=&\mathbb{E} \biggl( \exp-\lambda\sum
_{i\in\mathcal{I}; t_{i}\leq t} X^{i}_{t-t_{i}} \biggr)
= \exp \biggl(-\int_0^t \mathbb{N}
\bigl(1-\mathrm{e}^{-\lambda X_{t-s}} \bigr)\,\mathrm{d}s \biggr) %\mathrm{e}^{-x\imf{u_t}{\lambda}}
% \imf{\mathbb{P}_y}{df} \imf{\nu}{dy} \,\mathrm{d}s-\beta\int_0^t\int
\\
&=& %\mathrm{e}^{-x\imf{u_t}{\lambda}}
\exp \biggl(-\int_0^t\int
_0^\infty \bigl( 1-\mathrm{e}^{-yv_{t-s} ( \lambda ) } \bigr) \nu (\mathrm{d}y) \,\mathrm{d}s-\beta\int_0^t \int \bigl(
1-\mathrm{e}^{-\lambda x} \bigr) \eta_{t-s} (\mathrm{d}x) \,\mathrm{d}s \biggr)
\\
&=& %\mathrm{e}^ {-x\imf{u_t}{\lambda}}
\exp \biggl(-\int_0^t\Phi
\bigl( v_{t-s} ( \lambda) \bigr) \,\mathrm{d}s \biggr)=\exp \biggl(-\int
_0^t\Phi \bigl( v_{s} ( \lambda)
\bigr) \,\mathrm{d}s \biggr). %\\&=
\end{eqnarray*}

Hence, at least $Y$ has the correct one-dimensional distributions.
To prove that $Y$ is a CBI, we need to compute conditional
expectations. Let
\[
\mathcal{F}_t=\sigma \biggl\{\sum_{t_i\leq t}
\delta_{(
t_i,X^{i}_{t-t_i}) } \biggr\}.
\]
%
%Then\begin{eqnarray*}
%&=\mathbb{E}(\mathrm{e}^{-\lambda\sum_{t\leq s_n\leq t+s}\imf{f_n}{t+s-s_n} })
Then
\begin{eqnarray*}
\mathbb{E} \biggl[\exp \biggl(-\lambda\sum_{t\leq t_i\leq
t+s}X^{i}_{t+s-t_i}
\biggr)\Big\vert \mathcal{F}_t \biggr]&=&\mathbb{E} \biggl[\exp \biggl(-
\lambda\sum_{t\leq t_i\leq
t+s}X^{i}_{t+s-t_i}
\biggr) \biggr]
\\
&=&\exp \biggl(-\int_t^{t+s}\Phi \bigl(
v_{t+s-r} ( \lambda) \bigr) \,\mathrm{d}r \biggr)
\\
&=&\exp \biggl(-\int_{0}^{s}\Phi \bigl(
v_u ( \lambda) \bigr) \,\mathrm{d}u \biggr).
\end{eqnarray*}

On the other hand:
\begin{eqnarray*}
\mathbb{E} \biggl[\exp \biggl(-\lambda\sum_{t_i\leq t}
X^{i}_{t+s-t_i} \biggr)\Big\vert \mathcal{F}_t \biggr]
=\exp \biggl(-\sum_{t_{i}\leq t} X^{i}_{t-t_i}v_s(
\lambda) \biggr)
=\exp \bigl(-v_s(\lambda)Y_t \bigr).
\end{eqnarray*}
Hence,
\[
\mathbb{E} \bigl(\mathrm{e}^{-\lambda Y_{t+s}}|\mathcal{F}_t \bigr)=\exp
\biggl(-Y_t v_s(\lambda)-\int_{0}^{s}
\Phi \bigl(v_u(\lambda) \bigr) \,\mathrm{d}u \biggr).
\]
Because $Y$ is adapted to the filtration $ (\mathcal{F}_t,t\geq
0 )$ and since the conditional expectations of functionals of
$Y_{t+s}$ given $\mathcal{F}_t$ depend only on $Y_t$ we see that $Y$ is
a homogeneous Markov process. Since the Laplace transform of the
semigroup of $Y$ coincides with that of a
%We have thus proved that $Y$ is a homogeneous Markov process with
%respect to the filtration $\mathcal{F}_t,t\geq0$, to which $Y$ is
%adapted, and with the same semigroup as a
$\cbi(\Psi,\Phi)$, we conclude that $Y$ is a $\cbi(\Psi,\Phi)$ process
started at $0$.
\end{pf*}

%re4.1 #&#
\begin{rem}\label{grey2}
We mention that for the excursion measure $Q$ to exist (and then for
this spine decomposition to hold), one only needs that $\sigma>0$ or
$\int_{0}^{\infty}(1\wedge u) \nu_{1}(\mathrm{d}u)=\infty$. We refer the reader
to Duquesne and Labb\'e \cite{DuqLab}. Moreover, this spine
decomposition has also been used in the framework of stationary CBIs by
Bi \cite{Bi} to study the time of the most recent common ancestor.
%The Grey's condition fulfilled here by the map $\Psi$ is not necessary
%for this spine decomposition to hold. We refer the reader to Theorem
%1.1 of \cite{spinedecomposition}.
\end{rem}
\begin{pf*}{Proof of Theorem \ref{zeroset}}
From this spine decomposition, we establish a useful connection between
the random covering procedure and the set of the zeros of the CBI.
Results due to Fitzsimmons \textit{et al.} \cite{MR799145}, recalled in
Section~\ref{cutout}, will allow us to study the set $\mathcal{Z}$.
Denote for each $i\in\mathcal{I}$, $\zeta_{i}=\inf\{t\geq0;
X^{i}_{t}=0\}$. By standard properties of random Poisson measures, the
random measure $\sum_{i\in I}\delta_{(t_{i}, \zeta_{i})}$ is a Poisson
random measure with intensity $\mathrm{d}t\otimes\mathbb{N}(\zeta\in \mathrm{d}t)$.
Plainly, we have
\[
\{t\geq0, Y_{t}=0\}=\mathbb{R}_{+}{}\Big\backslash{} \biggl[ \bigcup
_{i\in
\mathcal{I}; X^i_0>0}[t_{i}, t_{i}+
\zeta_{i} [\, \cup\ \bigcup_{i\in\mathcal{I}; X^i_0=0}
]t_{i}, t_{i}+\zeta_{i}[ \biggr].
\]
The following crucial equality holds
%
%e3 #&#
%
\begin{equation}
\label{Z} \mathcal{Z}=\mathbb{R}_{+}{}\Big\backslash{}\bigcup_{i\in I}\,]t_{i},
t_{i}+\zeta_{i}[.
\end{equation}
Notice that some $t_{i}$ may belong to $\mathbb{R}_{+}\setminus\bigcup
_{i\in I}]t_{i}, t_{i}+\zeta_{i}[$ and not to $\{t\geq0, Y_{t}=0\}$.
However such point $t_{i}$ is a left accumulation point of the set $\{
t\geq0, Y_{t}=0\}$. Namely, consider $t_{i}\in\mathcal{Z}$, assume by
contradiction that there exists $\varepsilon>0$ such that
$[t_{i}-\varepsilon
, t_{i}[\,\cap\,\{t\geq0; Y_{t}=0\}=\varnothing$, then taking the closure we
have $[t_{i}-\varepsilon, t_{i}]\cap\mathcal{Z}=\varnothing$ which is
impossible since $t_{i}\in\mathcal{Z}$. Hence, equality \eqref{Z}
holds almost surely. Moreover, we have
\begin{eqnarray*}
\mathbb{N}(\zeta>t)&=&\beta Q[\zeta>t]+\int_{0}^{\infty}
\mathbb {P}_{x}(\zeta>t)\nu_{0}(\mathrm{d}x)
\\
&=&\beta v_{t}+\int_{0}^{\infty} \bigl(1-
\exp(-xv_{t}) \bigr)\nu _{0}(\mathrm{d}x)
\\
&=&\Phi(v_{t}).
\end{eqnarray*}
In order to establish the statement, we can directly use the results of
Section~\ref{cutout}. Part (i) and statement (a) of (ii) readily follow
from Theorem \ref{recovering}, taking
\[
\bar{\mu}(t)=\mathbb{N}(\zeta>t)=\Phi(v_{t})
\]
for all $t\geq0$. We get $\mathcal{Z}=\{0\}$ if and only if
\[
\int_{0}^{1}\exp \biggl[-\int
_{1}^{u}\Phi(v_{s})\,\mathrm{d}s \biggr]\,\mathrm{d}u=\infty.
\]
Using the differential equation satisfied by $v$, the substitution
$z=v_{t}$ and the condition $v_{0+}=\infty$ we get
\[
\int_{0}^{1}\exp \biggl[-\int
_{1}^{u}\Phi(v_{s})\,\mathrm{d}s \biggr]\,\mathrm{d}u=\int
_{v_{1}}^{\infty}\exp \biggl[\int_{v_{1}}^{z}
\frac{\Phi(u)}{\Psi
(u)}\,\mathrm{d}u \biggr]\frac{1}{\Psi(z)}\,\mathrm{d}z=\infty.
\]
Statement (b) follows from Proposition 1 of \cite{MR799145} (see also
Corollary 7.3 of \cite{subordinatorbertoin}). Namely the set $\mathcal
{Z}$ has a positive Lebesgue measure if and only if
\[
\int_{0}^{\infty}(s\wedge1)\mu(\mathrm{d}s)=\int
_{0}^{1}\Phi (v_{s})\,\mathrm{d}s<\infty.
\]
By the same substitution, this yields
\[
\int_{v_{1}}^{\infty}\frac{\Phi(s)}{\Psi(s)}\,\mathrm{d}s<\infty.
\]
The last statement (d) follows from Corollary 2 of \cite{MR799145}.
Namely, the necessary and sufficient condition for the uncovered set to
be a union of intervals is that the intensity measure $\mathbb
{N}(\zeta
\in \mathrm{d}t)$ has a finite mass. Therefore, assume that
$\Phi(v_{t})\longrightarrow_
{t \rightarrow0} c<\infty$. We have $v_{t}\longrightarrow_{t
\rightarrow0} \infty$ and
\[
\Phi(v_{t})=\beta v_{t}+\int_{0}^{\infty
}
\bigl(1-\mathrm{e}^{-v_{t}x} \bigr)\nu_{0}(\mathrm{d}x).
\]
By monotone convergence
\[
\int_{0}^{\infty} \bigl(1-\mathrm{e}^{-v_{t}x} \bigr)
\nu_{0}(\mathrm{d}x)\,\mathop{\longrightarrow}_{t
\rightarrow0}\, \nu_{0}
\bigl([0,\infty[ \bigr).
\]
Thus we get the conditions $\beta=0$ and $\nu_{0}([0,\infty[)<\infty$.
\end{pf*}
We now provide the proof of Lemma \ref{Hausdorff}. The arguments are
directly those used for the fractal dimensions of random cutout set.
\begin{pf*}{Proof of Lemma \ref{Hausdorff}} Combining Theorem 5.1 and
Corollary 7.6 of \cite{subordinatorbertoin}, we get that for every $t>0$,
\[
\overline{\operatorname{dim}} \bigl(\mathcal{Z}\cap[0,t] \bigr)=1-\mathop{\liminf
}_{u\rightarrow
0} \ \frac{\int_{u}^{1}\Phi(v_{s})\,\mathrm{d}s}{\log({1}/{u})} \qquad\mbox {a.s.}
\]
and
\[
\underline{\operatorname{dim}} \bigl(\mathcal{Z}\cap[0,t] \bigr)=1-\mathop{\limsup
}_{u\rightarrow
0} \ \frac{\int_{u}^{1}\Phi(v_{s})\,\mathrm{d}s}{\log({1}/{u})} \qquad\mbox {a.s.}
\]
We get the statement by substitution $u=v_{s}$. The second statement
concerning the largest zero is a direct application of Corollary 7.4 of
\cite{subordinatorbertoin}.
\end{pf*}
%
%s5 #&#
\section{Regularly varying branching-immigrating mechanisms}\label{RV}
In this section, we shall focus on specific mechanisms $\Psi$ and
$\Phi
$ such that the function \textit{ratio} $R(x)=\frac{\Phi(x)}{\Psi(x)}$
is regularly varying at $\infty$. The index of $R$, called in the
sequel, $\rho$ may be interpreted as representing the strength of the
immigration over the reproduction. For sake of conciseness, we only
work with a critical branching mechanism $\Psi$. A remarkable
phenomenon occurs when $\rho=-1$. In such case, other quantities are
involved. Surprisingly the quantities $\underline{r}:=
{\liminf}_{s\rightarrow\infty}\ sR(s)$ and $\overline{r}:=
{\limsup}_{s\rightarrow\infty}\ sR(s)$ play a crucial role.

We recall that a positive measurable function $f$ defined on some
neighbourhood of $\infty$ (resp.,~$0$) is said to be regularly
varying at $\infty$ (resp.,~$0$) with index $\rho$ if
for all $c>0$,
\[
f(cx)/f(x) \longrightarrow c^{\rho} \qquad\mbox{when }
x\rightarrow
\infty\mbox{ (resp., 0)}.
\]
%
%s5.1 #&#
\subsection{Polarity of zero and Hausdorff dimension of the zero set}
We first recall the definition of the upper and lower indices at
infinity of a Laplace exponent (we refer the reader to the seminal work
of Blumenthal and Getoor \cite{MR0123362} and to Duquesne and Le Gall's
work \cite{DuqLG}, page 557).
\begin{eqnarray*}
\underline{\operatorname{Ind}}(\Psi)&=&\sup \bigl\{ \alpha\geq0 ; \mathop{\lim
}_{\lambda
\rightarrow\infty} \Psi(\lambda)\lambda^{-\alpha}=\infty \bigr\} =\mathop{
\liminf}_{\lambda\rightarrow\infty}\frac{\log(\Psi(\lambda
))}{\log(\lambda)},
\\
\overline{\operatorname{Ind}}(\Psi)&=&\inf \bigl\{ \alpha\geq0 ; \mathop{\lim
}_{\lambda
\rightarrow\infty} \Psi(\lambda)\lambda^{-\alpha}=0 \bigr\}=\mathop{
\limsup}_
{\lambda\rightarrow\infty}\frac{\log(\Psi(\lambda))}{\log
(\lambda)}.
\end{eqnarray*}
This definition holds also for the Laplace exponent of a subordinator
$\Phi$, replace only $\Psi$ by $\Phi$ above.
%then $\rho=\underline{\mbox{Ind}}(\Phi)-\overline{\mbox{Ind}}(\Psi)$.
%Furthermore, we have $\underline{\mbox{Ind}}(\Phi)-\overline{
%is a slowly varying function $l$ such that $R(x)=x^{\rho}l(x)$.
%Moreover we have that $\frac{\log(l(x))}{\log(x)} {x
%1.3.6-(i) of \cite{regularvariation}), therefore
%$$\frac{\log(R(x))}{\log(x)}=\rho+\frac{\log(l(x))}{\log(x)}\underset{x
%Plainly $$\log(R(x))/\log(x)=\left(\log(\Phi(x))-\log(\Psi(x))
%$$\underset{x\rightarrow\infty}\lim\frac{\log(R(x))}{\log(x)}=
%We deduce then that
%$$\log\left(\frac{1}{x}\right)-\log\left(\frac{1}{x}\right)\frac{
%Since $\log\left(\frac{l(x)}{l(\lambda x)}\right)\underset{x
%the limit of the quantity in the right-hand side of the last equality
%is the same as that of
%$$\frac{\rho\log(\lambda)\log(x)}{\rho\log(x)+\log(l(x))}=\frac{\log(
%Since the function $l$ is slowly varying, we have that $\frac{
%0$ . Finally,
%since\[
%$\rho=\overline{\mbox{Ind}}(\Phi)-\underline{\mbox{Ind}}(\Psi)=
%
%th9 #&#
\begin{thmm} \label{polarityRV}
Assume that $R$ is regularly varying at $+\infty$ with index $\rho$.
\begin{longlist}[(iii)]
\item[(i)] If $\rho>-1$, then $0$ is polar.
\item[(ii)] If $\rho<-1$, then $0$ is not polar and the zero set is heavy.
\item[(iii)] If $\rho=-1$. Define the quantities
\[
\overline{r}:=\mathop{\limsup}_{s\rightarrow\infty}\ s\frac{\Phi
(s)}{\Psi(s)}
\]
and
\[
\underline{r}:=\mathop{\liminf}_{s\rightarrow\infty}\ s\frac{\Phi
(s)}{\Psi(s)}.
\]
Both quantities belong to $[0,\infty]$.
\begin{longlist}[(a)]
\item[(a)]If $\overline{r}<\underline{\operatorname{Ind}}(\Psi)-1$,
then $0$ is
not polar. Moreover if $\underline{r}>0$, then the zero set is light.
\item[(b)]If $\underline{r}\geq\overline{\operatorname{Ind}}(\Psi
)-1$, then $0$
is polar.
\end{longlist}
\end{longlist}
\end{thmm}
\begin{pf} Assume that $\rho>-1$, therefore $z\mapsto\int_{\theta
}^{z}R(u)\,\mathrm{d}u$ is regularly varying with index $\rho+1$ and by Karamata's
theorem (see, e.g., Proposition 1.5.8 of \cite
{regularvariation}) we have
%
%e4 #&#
%
\begin{equation}
\label{prim} \mbox{For all } z\geq\theta,\qquad \int_{\theta}^{z}R(u)\,\mathrm{d}u=C+z^{\rho+1}l(z)
\end{equation}
with $C$ a constant and $l$ a slowly varying function. By Corollary
\ref
{recurrence}, we have to study the following integral
\[
\mathcal{I}:=\int_{\theta}^{\infty}\exp \biggl(\int
_{\theta
}^{z}R(u)\,\mathrm{d}u \biggr)\frac{1}{\Psi(z)}\,\mathrm{d}z.
\]

\begin{enumerate}[(iii)]
\item[(i)] If $\rho>-1$, since $l$ is slowly varying, we have
$z^{(\rho
+1)/2}l(z)\longrightarrow_{z\rightarrow\infty} +\infty$,
moreover observing that $\Psi(z)=\mathrm{O}(z^{2})$, we have clearly that
$\mathcal{I}=\infty$.
\item[(ii)]If $\rho<-1$, by Proposition 1.5.10 of \cite
{regularvariation}, we have $\int_{\theta}^{\infty}R(u)\,\mathrm{d}u<\infty$, then
$\mathcal{I}<\infty$ and furthermore by statement (ii)(b) of Theorem
\ref
{zeroset}, the zero set has a positive Lebesgue measure. %To obtain the
%complete statement $(ii)$, it remains to study the following integral
%$$\int_{0}^{\theta}\exp\left[-\int_{x}^{\theta}R(u)\,\mathrm{d}u\right]\frac{dx}{
%For $x$ small enough, we have $e^{x^{\rho+1}l(x)}/\Psi(x)\geq1/
%$(iii)$ of Corollary \eqref{recurrence} is fulfilled.
%
\item[(iii)] In the case $\rho=-1$, the map $R$ is regularly varying
with index $-1$. We mention that integrals of such regularly varying
functions yields to de Haan functions (see Chapter~3 of \cite
{regularvariation}).
\begin{enumerate}
\item[(a)] Let $\varepsilon>0$. For $s$ large enough, $R(s)\leq\frac
{\overline{r}}{s}+\frac{\varepsilon}{s}$. Therefore for large enough $z$,
\[
\exp \biggl( \int_{\theta}^{z}R(s)\,\mathrm{d}s \biggr)
\frac{1}{\Psi(z)}\leq C\frac
{z^{\overline{r}+\varepsilon}} {\Psi(z)}.
\]
Moreover, $\Psi(z)\geq C'z^{\underline{\operatorname{Ind}}(\Psi
)-\varepsilon}$ and thus
\[
\exp \biggl( \int_{\theta}^{z}R(s)\,\mathrm{d}s \biggr)
\frac{1}{\Psi(z)}\leq C''z^{\overline{r}-\underline{\operatorname{Ind}}(\Psi)+2\varepsilon}.
\]
As $\varepsilon$ is arbitrarily small, if $\overline{r}-\underline
{\operatorname
{Ind}}(\Psi)<-1$ then
\[
\int_{\theta}^{\infty}z^{\overline{r}-\underline{\operatorname
{Ind}}(\Psi
)+2\varepsilon}\,\mathrm{d}z<\infty.
\]
This implies that $0$ is not polar. Assume $\underline{r}>0$ for $s$
large enough, $R(s)\geq\underline{r}/2s$, we have clearly $\int_{\theta
}^{\infty}R(s)\,\mathrm{d}s=\infty$.
\item[(b)] For $s$ large enough, $sR(s)\geq\underline{r}-\varepsilon$ and
then $\exp ( \int_{\theta}^{z}R(s)\,\mathrm{d}s )\frac{1}{\Psi
(z)}\geq C
\frac{z^{\underline{r}-\varepsilon}}{\Psi(z)}$
By the same reasoning, we have $\Psi(z)\leq C'z^{\overline
{\operatorname
{Ind}}(\Psi)+\varepsilon}$ and we get
\[
\exp \biggl( \int_{\theta}^{z}R(s)\,\mathrm{d}s \biggr)
\frac{1}{\Psi(z)}\geq C'' z^{\underline{r}-\overline{\operatorname{Ind}}(\Psi)-2\varepsilon}.
\]
If $\underline{r}-\overline{\operatorname{Ind}}(\Psi)>-1$, since
$\varepsilon$ is
arbitrarily small, then $\int_{\theta}^{\infty}z^{\underline
{r}-\overline{\operatorname{Ind}}(\Psi)-2\varepsilon}\,\mathrm{d}z=\infty$ and
$\mathcal
{I}=\infty$.

If $\underline{r}-\overline{\operatorname{Ind}}(\Psi)=-1$, the direct
computation
\[
\int_{\theta}^{\infty}z^{-1-2\varepsilon}\,\mathrm{d}z=\frac{\theta^{-2\varepsilon
}}{2\varepsilon}
\]
yields the lower bound $\mathcal{I}\geq C''\frac{\theta^{-2\varepsilon
}}{2\varepsilon}$. Letting $\varepsilon$ going to $0$, we get $\mathcal
{I}=\infty$. \qed%\qedhere
\end{enumerate}
\end{enumerate}
\noqed
\end{pf}
%
%In some cases, provided that the behavior of $R$ in the neighbourhood
%of $\infty$ is \textit{simple}, we are able to provide the Hausdorff
%dimension of the zero set.
%denote by $r$ this quantity. If zero is not polar then for all $t>0$,
%$$\mbox{dim}_{H}(\mathcal{Z}\cap[0,t])=1-\frac{r}{\overline{
%
%pr10 #&#
\begin{prop} \label{dim} Assume that $\overline{r}<\underline
{\operatorname
{Ind}}(\Psi)-1$ then for all $t>0$,
\[
1-\frac{\overline{r}}{\overline{\operatorname{Ind}}(\Psi)-1}\leq \overline{\operatorname {dim}} \bigl(\mathcal{Z}
\cap[0,t] \bigr)\leq1-\frac{\underline{r}}{\overline
{\operatorname
{Ind}}(\Psi)-1}
\]
and
\[
1-\frac{\overline{r}}{\underline{\operatorname{Ind}}(\Psi)-1}\leq \underline {\operatorname{dim}} \bigl(\mathcal{Z}
\cap[0,t] \bigr)\leq1-\frac{\underline
{r}}{\underline{\operatorname{Ind}}(\Psi)-1}.
\]
\end{prop}
Before tackling the proof, we need a lemma.
%
%le11 #&#
\begin{lem}[(Duquesne Le Gall, Lemma 5.6 of \cite{DuqLG})] The following
equalities hold
\[
\mathop{\liminf}_{t \rightarrow0}\frac{\log(v_{t})}{\log(1/t)}=\frac
{1}{\overline{\operatorname{Ind}}(\Psi)-1}
\]
and
\[
\mathop{\limsup}_{t \rightarrow0}\frac{\log(v_{t})}{\log(1/t)}=\frac
{1}{\underline{\operatorname{Ind}}(\Psi)-1}.
\]
\end{lem}
\begin{pf} Duquesne and Le Gall established
\[
\mathop{\liminf}_{t \rightarrow0}\frac{\log(v_{t})}{\log(1/t)}\leq \frac{1}{\overline{\operatorname{Ind}}(\Psi)-1}
\]
and
\[
\mathop{\limsup}_{t \rightarrow0}\frac{\log(v_{t})}{\log(1/t)}\leq \frac{1}{\underline{\operatorname{Ind}}(\Psi)-1}.
\]
The other inequalities are also true. We provide here a proof. Denote
$\eta:=\overline{\operatorname{Ind}}(\Psi)$, there exists $\varepsilon>0$
arbitrarily small such that
\[
q^{-(\eta+\varepsilon)}\Psi(q)\,\mathop{\longrightarrow}_{q\rightarrow\infty
}\, 0.
\]
For all $C>0$, for large enough $q$ we have $\frac{1}{\Psi(q)}\geq
\frac
{C}{q^{\eta+\varepsilon}}$. Recall that $\eta>1$, therefore for small
enough $t$, we have
\[
\int_{v_{t}}^{\infty}\frac{\mathrm{d}q}{\Psi(q)}=t\geq C\int
_{v_{t}}^{\infty
}\frac{\mathrm{d}q}{q^{\eta+\varepsilon}}=\frac{C}{\eta+\varepsilon
-1}v_{t}^{1-\eta
-\varepsilon}.
\]
This yields
\[
\log \biggl(\frac{1}{t} \biggr)\leq\log \biggl(\frac{\eta+\varepsilon
-1}{C} \biggr)
+(\eta-1+\varepsilon)\log(v_{t}).
\]
We then have
\[
\mathop{\liminf} _{t\rightarrow0}\frac{\log(v_{t})}{\log(1/t)}\geq \frac{1}{\eta-1+\varepsilon}.
\]
As $\varepsilon$ is arbitrarily small, the result follows and the first
equality in the statement is established.

We study now the second equality. Denote $\gamma:=\underline
{\operatorname
{Ind}}(\Psi)$. Let $\gamma'>\gamma$. There exists a sequence $(u_{n},
n\geq1)$ such that $u_{n}\,{\longrightarrow}_{n\rightarrow\infty}\,
\infty$ and $\Psi(2u_{n})\leq2^{\gamma'}
u_{n}^{\gamma'}$. For $u\in\,]0, 2u_{n}]$, we have by convexity $\frac
{\Psi(u)}{u}\leq\frac{\Psi(2u_{n})}{2u_{n}}$ and therefore $\Psi
(u)\leq2^{\gamma'-1}uu_{n}^{\gamma'-1}$ when $u\in\,]0,2u_{n}]$. Define
the function $F(a)=\int_{a}^{\infty}\frac{\mathrm{d}u}{\Psi(u)}$, we have
\[
F(u_{n})=\int_{u_{n}}^{\infty}
\frac{\mathrm{d}u}{\Psi(u)}\geq\int_{u_{n}}^{2u_{n}}
\frac{\mathrm{d}u}{\Psi(u)}\geq2^{1-\gamma'}u_{n}^{1-\gamma
'}\log(2).
\]
Then
\[
\mathop{\liminf}_{a\rightarrow\infty} \frac{\log(1/F(a))}{\log
(a)}\leq\mathop{\lim}_{n\rightarrow\infty}
\frac{\log(1/F(u_{n}))} {
\log(u_{n})}\leq\gamma'-1.
\]
As on page 592 of \cite{DuqLG}, observe that by definition of $v_{t}$:
\[
\biggl(\mathop{\limsup}_{t\rightarrow0} \frac{\log(v_{t})}{\log
(1/t)} \biggr)^{-1}=\mathop{
\liminf}_{a \rightarrow\infty}\frac{\log
(1/F(a))}{\log(a)}.
\]
We deduce that
\[
\mathop{\limsup}_{t\rightarrow0} \frac{\log(v_{t})}{\log(1/t)}\geq \frac{1}{\gamma'-1},
\]
by letting $\gamma'$ go to $\gamma$, we obtain the wished inequality.
\end{pf}
\begin{pf*}{Proof of Proposition \ref{dim}} Let $\varepsilon>0$, by
assumption, we have for $C$ large enough
\[
\mathop{\sup}_{s\in[C,\infty[} sR(s) \leq\overline{r}+\varepsilon \quad\mbox{and}\quad
\mathop{\inf}_{s\in[C,\infty[} sR(s) \geq\underline{r}-\varepsilon.
\]
Therefore,
\[
(\underline{r}-\varepsilon) \bigl[\log(v_{t})-\log(C) \bigr]\leq\int
_{C}^{v_{t}}\frac
{\Phi(u)}{\Psi(u)}\,\mathrm{d}u \leq(\overline{r}+
\varepsilon) \bigl[\log(v_{t})-\log(C) \bigr].
\]
By using the previous lemma and Lemma \ref{Hausdorff}, we plainly get
\[
1-\frac{\overline{r}+\varepsilon}{\overline{\operatorname{Ind}}(\Psi
)-1}\leq \overline{\operatorname{dim}} \bigl(\mathcal{Z}
\cap[0,t] \bigr)\leq1-\frac
{\underline
{r}-\varepsilon}{\overline{\operatorname{Ind}}(\Psi)-1}.
\]
As $\varepsilon$ is arbitrarily small, we get the statement. Same
arguments hold for the lower box-counting dimension.
\end{pf*}

%s5.2 #&#
\subsection{Recurrence and regular variation at $0+$}
To study the recurrence of zero, we need information on the behaviour
of the map $R$ in the neighbourhood of $0+$ (see statement (iii) in
Corollary \ref{recurrence}). In the same vein as our previous result on
polarity, a natural assumption is to consider the map $R$ with regular
variation at $0+$.

In order to state the following result, we need to introduce the lower
and upper indices of a Laplace exponent at $0+$:
\begin{eqnarray*}
\underline{\operatorname{ind}}(\Psi)&=&\sup \bigl\{ \alpha\geq0 ; \mathop{\lim
}_{\lambda
\rightarrow0} \Psi(\lambda)\lambda^{-\alpha}=0 \bigr\} = \mathop{
\liminf}_{\lambda\rightarrow0}\frac{\log(\Psi(\lambda))}{\log
(\lambda
)},
\\
\overline{\operatorname{ind}}(\Psi)&=&\inf \bigl\{ \alpha\geq0 ; \mathop{\lim
}_{\lambda
\rightarrow0} \Psi(\lambda)\lambda^{-\alpha}=\infty \bigr\} =\mathop{
\limsup}_{\lambda\rightarrow0}\frac{\log(\Psi(\lambda))}{\log
(\lambda)}.
\end{eqnarray*}

%We begin by stating an easy result concerning these quantities. As
%previously, we shall consider only conservative processes.
%observe that $1\leq\underline{\mbox{ind}}(\Psi)\leq\overline{
%varying at $0+$ then $\overline{\mbox{ind}}(\Psi)=\underline{
%(however this example does not satisfy Grey's condition).
%$\underline{\mbox{ind}}(\Psi)\geq1$. Denote $\beta=\underline{
%assume that $x^{-\beta-\varepsilon}\Psi(x)\underset{x\rightarrow0}{
%quantity $\varepsilon$ being arbitrary, $\beta+\varepsilon<1$ corresponds
%with $\beta<1$ and this implies plainly that $\int_{0+}\frac{dx}{
%
%th12 #&#
\begin{thmm}\label{recurrenceRV} Assume that the map $R$ is regularly
varying at $+\infty$ with index $\rho$ and at $0+$ with index $\kappa$.
If $\rho<-1$ or $\rho=-1$ and $\overline{r}<\underline
{\operatorname{Ind}}(\Psi
)-1$ (then $0$ is not polar) and
\begin{enumerate}[(iii)]
\item[(i)] if $\kappa<-1$ then $0$ is transient,
\item[(ii)]if $\kappa>-1$ then $0$ is recurrent,
\item[(iii)]if $\kappa=-1$. Define the quantities
\[
\overline{\kappa}:=\mathop{\limsup}_{x\rightarrow0}\ xR(x)
\]
and
\[
\underline{\kappa}:=\mathop{\liminf}_{x\rightarrow0}\ xR(x).
\]
Both quantities belong to $[0,\infty]$.
\begin{enumerate}[(a)]
\item[(a)] If $\overline{\kappa}-\underline{\operatorname
{ind}}(\Psi)\leq-1$,
then $0$ is recurrent.
\item[(b)] If $\underline{\kappa}-\overline{\operatorname
{ind}}(\Psi)> -1$, then
$0$ is transient.
\end{enumerate}
\end{enumerate}
\end{thmm}
\begin{pf} We assume that $0$ is not polar, then we have to study
the following integral
\[
\mathcal{J}:=\int_{0}^{\theta}\exp \biggl[-\int
_{x}^{\theta
}R(u)\,\mathrm{d}u \biggr]\frac{\mathrm{d}x}{\Psi(x)}.
\]
Assume that $R$ is regularly varying at $0$ with index $\kappa\in
\mathbb{R}$, then by an easy adaptation of Propositions 1.5.8 and
1.5.10 of \cite{regularvariation} to the setting of regular variation
at $0+$, we get
\[
\mathop{\lim}_{x\rightarrow0} \int_{x}^{\theta
}R(u)\,\mathrm{d}u= \cases{
=\infty,&\quad$ \mbox{if } \kappa<-1,$\vspace*{2pt}\cr
<\infty,&\quad$\mbox{if } \kappa>-1.$}
\]
More precisely, if $\kappa<-1$, we have $\int_{x}^{\theta
}R(u)\,\mathrm{d}u=x^{\kappa+1}l(x)$ with $l$ a slowly varying function at $0+$.
\begin{itemize}
\item Assume $\kappa<-1$. Let $\varepsilon>0$, we have for small enough
$x$, $\Psi(x)\geq Cx^{\overline{\operatorname{ind}}(\Psi)+\varepsilon
}$. We shall
prove that the map
\[
\label{quantity} x \mapsto\frac{\exp(-x^{\kappa
+1}l(x))}{x^{\overline{\operatorname{ind}}(\Psi)+\varepsilon}}
\]
is bounded in the neighbourhood of $0$. This will imply that $\mathcal
{J}$ is finite. Taking the logarithm, we have
\[
\log \biggl(\frac{\exp(-x^{\kappa+1}l(x))}{x^{\overline
{\operatorname{ind}}(\Psi
)+\varepsilon}} \biggr)=-x^{\kappa+1}l(x) \biggl(1+ \bigl(
\overline {\operatorname{ind}}(\Psi )+\varepsilon \bigr)\frac{1}{l(x)}
\frac{\log(x)}{x^{\kappa+1}} \biggr).
\]
On the one hand, $x^{\kappa+1}l(x)\longrightarrow_{x \rightarrow
0}+\infty$, on the other hand the map\vspace*{1pt} $x\mapsto\frac
{\log(x)}{l(x)}$ is slowly varying at $0+$, and thus using Potter's
bound (see, e.g., Theorem 1.5.6(iii) in \cite{regularvariation}), we
have $x^{-(\kappa+1)}\frac{\log(x)}{l(x)} \,{\longrightarrow}_{x
\rightarrow
0}\, 0$. Finally,
\[
-x^{\kappa+1}l(x) \biggl(1+ \bigl(\overline{\operatorname{ind}}(\Psi )+
\varepsilon \bigr)\frac
{1}{l(x)}\frac{\log(x)}{x^{\kappa+1}} \biggr)\mathop{\longrightarrow
}_{x\rightarrow
0} -\infty
\]
and the map is bounded. We deduce that $\mathcal{J}<\infty$ and then
$0$ is transient.
\item If $\kappa>-1$, provided that $\int_{0+}\frac{\mathrm{d}x}{\Psi
(x)}=\infty
$, we have $\mathcal{J}=\infty$ and $0$ is recurrent.
\item We deal now with the case $\kappa=-1$. We prove first the
recurrence criterion (a). Let $\varepsilon>0$, for $u$ small enough, we
have
\[
R(u)\leq\frac{\overline{\kappa}}{u}+\frac{\varepsilon
}{u} \quad\mbox{and}\quad \Psi(u) \leq Cu^{\underline{\operatorname{ind}}(\Psi
)-\varepsilon}.
\]
We deduce that
\[
\frac{1}{\Psi(x)}\exp \biggl(-\int_{x}^{\theta}R(u)\,\mathrm{d}u
\biggr) \geq Cx^{\overline{\kappa}-\underline{\operatorname{ind}}(\Psi
)+2\varepsilon}.
\]
If $\overline{\kappa}-\underline{\operatorname{ind}}(\Psi)<-1$,
since $\varepsilon
$ is arbitrarily small, we have $\mathcal{J}=\infty$.

If $\overline{\kappa}-\underline{\operatorname{ind}}(\Psi)=-1$, then
\[
\int_{0}^{\theta}x^{-1+2\varepsilon}\,\mathrm{d}x=
\frac{\theta^{2\varepsilon
}}{2\varepsilon},
\]
and $\mathcal{J}\geq C\frac{\theta^{2\varepsilon}}{2\varepsilon}$ and letting
$\varepsilon$ going to $0$, we obtain $\mathcal{J}=\infty$. Therefore, $0$
is recurrent. We prove now statement (b). Assume $\underline{\kappa
}-\overline{\operatorname{ind}}(\Psi)> -1$. We have for $u$ small enough
$uR(u)\geq\underline{\kappa}-\varepsilon$ and $\Psi(u)\geq
u^{\overline
{\operatorname{ind}}(\Psi)+\varepsilon}$. Therefore,
\[
\frac{1}{\Psi(x)}\exp \biggl(-\int_{x}^{\theta}R(u)\,\mathrm{d}u
\biggr)\leq Cx^{\underline{\kappa}-\overline{\operatorname{ind}}(\Psi
)-2\varepsilon}.
\]
Since $\varepsilon$ is arbitrarily small, we can choose one such that
$\underline{\kappa}-\overline{\operatorname{ind}}(\Psi)-2\varepsilon
>-1$. This
implies $\mathcal{J}<\infty$ and the transience follows. \qed%\qedhere
\end{itemize}
\noqed
\end{pf}
%
%re5.1 #&#
\begin{rem} If $\kappa>-1$, then $\int_{0+}R(u)\,\mathrm{d}u<\infty$ and by
Theorem \ref{stationary}, the CBI has a stationary law.
\end{rem}
We study in the sequel the specific case of stable and gamma mechanisms.
%s5.2.1 #&#
\subsubsection{Stable and gamma mechanisms}
Consider $(\Psi, \Phi)$ of the form $\Psi(q)={d}q^{\alpha}$ and $\Phi
(q)=d'q^{\beta}$ with $\alpha\in(1,2]$, $\beta\in(0,1]$ and $d,
d'\in
(0,\infty)$. The CBI process $(Y_{t}, t\geq0)$ associated is said to
be \textit{stable}. Obviously, we have $\overline{\operatorname
{Ind}}(\Psi
)=\alpha$, $\underline{\operatorname{Ind}}(\Phi)=\beta$. The map
$R$ is
regularly varying at $0+$ and at $+\infty$ with index $\rho=\kappa
=\beta
-\alpha$. We can thus apply the previous results in Theorem \ref
{polarityRV} and Theorem~\ref{recurrenceRV}.
\begin{itemize}
\item If $\beta>\alpha-1$, then $0$ is polar,
\item if $\beta<\alpha-1$, then $0$ is transient and the zero set is heavy,
\item if $\beta=\alpha-1$, we have $r:=\lim_{x\rightarrow\infty
}xR(x)=\frac{d'}{d}$ and $\kappa=-1$. Two cases may occur.
\begin{itemize}[--]
\item[--] If $\frac{d'}{d}\geq\alpha-1$, then $0$ is polar,
\item[--] if $\frac{d'}{d}< \alpha-1$, then $0$ is recurrent and by
Proposition \ref{dim}
\[
\operatorname{dim}_{H} \bigl(\mathcal{Z}\cap[0,t] \bigr)=1-
\frac{1}{\alpha
-1}\frac{d'}{d}.
\]
\end{itemize}
\end{itemize}
Notice that in this stable framework, we cannot have $\rho<-1$ and
$\kappa\geq-1$.
%
%pr13 #&#
\begin{prop}
Let $(Y_{t}, t\geq0)$ denote a stable critical CBI started at $0$ with
parameters $\alpha, \beta$ satisfying $\beta=\alpha-1$.
\[
\mathcal{Z}=\overline{\{t\geq0; Y_{t}=0\}}=\overline{\{
\sigma_{t}, t\geq0\}}
\]
with $(\sigma_{t}, t\geq0)$ a stable subordinator with index $\gamma
=1-\frac{1}{\alpha-1}\frac{d'}{d}$.
\end{prop}
\begin{pf}
In order to get that the subordinator $(\sigma_{t}, t\geq0)$ is a
$\gamma$-stable one, we shall use self-similarity result. The
self-similarity of $(Y_{t}, t\geq0)$ follows by inspection (see also
\cite{Patie}). Namely, we have
\[
\mathbb{E}_{0} \bigl[\mathrm{e}^{-qY_{t}} \bigr]=\exp \biggl(-\int
_{0}^{t}\Phi \bigl(v_{s}(q) \bigr)\,\mathrm{d}s
\biggr)
\]
with
$v_{s}(q)=q[1+d(\alpha-1)q^{\alpha-1}s]^{-{1}/{(\alpha-1)}}$. An easy
computation yields that the processes $(kY_{t}, t\geq0)$ and
$(Y_{k^{\alpha-1}t}, t\geq0)$ have the same law. We deduce that the
regenerative set $\mathcal{Z}$ is self-similar meaning that for any
$k>0$, $k\mathcal{Z}\stackrel{\mathrm{law}}{=}\mathcal{Z}$. The only regenerative
sets satisfying this property are the closure of the range of stable
subordinator (see, e.g., Section~3.1.1 of \cite
{subordinatorbertoin}). Proposition \ref{dim} provides the Hausdorff
dimension of the set and therefore the index of stability (see Theorem
5.1 of \cite{subordinatorbertoin}).
\end{pf}
%
%Let $g_{t}=\sup\{s<t; s\in\mathcal{Z}\}$ denotes the last passage
%time in $\mathcal{Z}$ before time $t>0$. By a direct application of
%Proposition 3.1 of \cite{subordinatorbertoin}, $g_{1}$ has the
%generalized arcsine law
%$$\mathbb{P}(g_{1}\in ds)=\frac{\sin(\gamma\pi)}{\pi}s^{

Consider the immigration mechanism $\Phi(q)=\frac{\Gamma(\beta
+q)}{\Gamma(\beta)\Gamma(q)}{\sim}_{q\rightarrow\infty}\frac
{q^{\beta}}{\Gamma(\beta)}$. The subordinator with Laplace exponent
$\Phi$ is called Lamperti stable subordinator (see, e.g.,
\cite
{MR2792485}). The map $\Phi$ is regularly varying at $+\infty$ with
index $\beta$ and at $0+$ with index $1$. Assume that the reproduction
mechanism $\Psi$ is $\alpha$-stable, as previously we can apply
Theorems \ref{polarityRV} and \ref{recurrenceRV}. We easily get
\begin{itemize}
\item if $\beta>\alpha-1$, then $0$ is polar,\vadjust{\goodbreak}
\item if $\beta<\alpha-1$, then $0$ is recurrent and the zero set is
heavy,
\item if $\beta=\alpha-1$, we have $r:={\lim}_{s\rightarrow\infty
}\,sR(s)=\frac{1}{{d}\Gamma(\alpha-1)} \mbox{ and } \kappa=1-\alpha
\geq-1$. Two cases may occur.
\begin{itemize}[--]
\item[--] if $d\leq\frac{1}{\Gamma(\alpha)}$ then $0$ is polar,
\item[--] if $d>\frac{1}{\Gamma(\alpha)}$ then $0$ is not polar, and we
easily verify that $0$ is recurrent:
\begin{eqnarray*}
\mbox{if } \alpha&=&2 \mbox{ then } \underline{\kappa}=\overline {\kappa }=1,
\mbox{ and } 0 \mbox{ is recurrent},
\\
\mbox{if } \alpha&\in&(1,2) \mbox{ then } \underline{\kappa }=\overline {
\kappa}=0, \mbox{ and } 0 \mbox{ is recurrent}.
\end{eqnarray*}
\end{itemize}
Finally, by Proposition \ref{dim} we get
\[
\operatorname{dim}_{H} \bigl(\mathcal{Z}\cap[0,t] \bigr)=1-
\frac{1}{\Gamma
(\alpha)d}.
\]
\end{itemize}
Last, consider now the case of a Gamma immigration mechanism and a
stable branching one. Let $a>0, b>0$ and $\alpha\in(1,2], d>0$.
\begin{eqnarray*}
\Phi(x)&=&a\log(1+x/b)=\int_{0}^{\infty
}
\bigl(1-\mathrm{e}^{-xu} \bigr)au^{-1}\mathrm{e}^{-bu}\,\mathrm{d}u,
\\
\Psi(x)&=&{d}x^{\alpha}.
\end{eqnarray*}
We can observe that $\Phi$ is slowly varying at $\infty$ and regularly
varying at $0+$ with index $1$. Thus $R\dvtx x\mapsto\frac{\Phi(x)}{\Psi
(x)}$ is regularly varying at $0$ with index $\kappa=1-\alpha$ and at
$+\infty$ with index $\rho=-\alpha$. Applying Theorems \ref{polarityRV}
and \ref{recurrenceRV}, we obtain that the zero set is heavy and further:
\begin{itemize}
\item if $\alpha\in(1,2)$, then $\kappa>-1$ and $0$ is recurrent,
\item if $\alpha=2$, then $\kappa=-1$ and $\overline{\kappa
}=\underline
{\kappa}=\frac{a}{bd}$. Therefore, $0$ is recurrent if $\frac
{a}{b}\leq
d$, otherwise $0$ is transient.
\end{itemize}
%
%s6 #&#
\section{Ornstein--Uhlenbeck processes}\label{OU}
We recall here basics on Markov processes of Ornstein--Uhlenbeck type.
Contrary to CBI processes, these processes in general take values in
$\mathbb{R}$. Consider $x\in\mathbb{R}$, $\gamma\in\mathbb{R}_{+}$,
and $(A_{t}, t\geq0)$ a one-dimensional L\'evy process with
characteristic function given by $\eta$ such that
\begin{eqnarray*}
\mathbb{E} \bigl[\mathrm{e}^{\mathrm{i}zA_{t}} \bigr]&=&\exp \bigl(t\eta(z) \bigr),
\\
\eta(z)&=&-\frac{\sigma^{2}}{2}z^{2}+\mathrm{i}bz+\int_{\mathbb
{R}}
\bigl(\mathrm{e}^{\mathrm{i}zx}-1-\mathrm{i}zx1_{\{|x|\leq1\}} \bigr)\nu(\mathrm{d}x),
\end{eqnarray*}
where $\sigma\geq0, b\in\mathbb{R}$ and $\int_{\mathbb
{R}}(|x|^{2}\wedge1) \nu(\mathrm{d}x)<\infty$.
A process $(X_{t}, t\geq0)$ is said to be an OU type process if it
satisfies the following equation
\[
X_{t}=x-\gamma\int_{0}^{t}X_{s}\,\mathrm{d}s+A_{t}.\vadjust{\goodbreak}
\]
Generalized Ornstein--Uhlenbeck processes valued in $\mathbb{R}_{+}$
belong actually to the class of CBI processes. When $A$ is subordinator
with L\'evy measure $\nu$ and drift $d$, the Ornstein--Uhlenbeck is a
CBI with $\Psi(z)=\gamma z$ and $\Phi(z)=\mathrm{d}z+\int_{0}^{\infty
}(1-\mathrm{e}^{-zx})\nu(\mathrm{d}x)$. In this case, Grey's condition is not fulfilled
and clearly the state $0$ is polar. Most of the Ornstein--Uhlenbeck
processes are not CBI processes (we refer, {e}.g., to the
discussion about their stationary laws in Proposition 4.7 in
Keller-Ressel and Mijatovi\'c \cite{Keller}).
% it is worth noticing that the zero sets of OU processes are
%infinitely divisible.
%Really?
%not at all
However, there is an interesting class of OU processes whose zero sets
are random cutout sets: Ornstein--Uhlenbeck processes whose driving L\'
evy process is self-similar of index $\alpha$. A L\'evy process $A$ is
self-similar of index $\alpha$ if for any $c,t>0$ there is equality in
law between $A_{ct}$ and $c^{1/\alpha}A_t$. Note that it excludes the
asymmetric Cauchy process $(\alpha=1)$ which is not strictly stable.
The corresponding characteristic function $\eta$ is rather involved and
we refer the reader to page 11 of Kyprianou \cite{MR2250061}.
%
%th14 #&#
\begin{thmm}\label{stableOU} Let $(X_{t}, t\geq0)$ be an
Ornstein--Uhlenbeck process started at $0$ driven by a self-similar L\'
evy process of index $\alpha$.
%Suppose $A$ is an $\alpha$-stable process.
If $\alpha\in(0,1]$, then the zero set of $X$ equals $\{0\}$ almost
surely. If $\alpha\in(1,2]$, then the zero set of $X$ is a random
cutout set whose cutting measure has density with respect to Lebesgue
measure given by
%$$
%(\tilde t,\tilde x)\mapsto(1-\beta)\mathrm{e}^{\tilde x}/(\mathrm{e}^{\tilde x}-1)^2 d
%$$where $\beta=1-1/\alpha$.
%
\[
z\mapsto(1-\beta)\mathrm{e}^{z}/ \bigl(\mathrm{e}^{z}-1 \bigr)^2,
\]
where $\beta=1-1/\alpha$.
\end{thmm}
%
%Note concerning Truman and Pitman-Yor.
%
\begin{pf}
Consider the process
\[
\tilde X_t= Ce^{-\gamma t}A_{\mathrm{e}^{\gamma t \alpha}},\qquad\mbox{where } C=
\frac{1}{( \alpha\gamma) ^{1/\alpha}}.
\]
Because $A$ is self-similar, it follows that $\tilde X$ is stationary.
Using integration by parts, it follows that
\[
\mathrm{d}\tilde X_t=Ce^{-\gamma t} \,\mathrm{d}A_{\mathrm{e}^{\gamma t \alpha}}-\gamma\tilde
X_t \,\mathrm{d}t.
\]
%
% Indeed, the covariation of $\mathrm{e}^{-\gamma\cdot}$ and $A_{\mathrm{e}^{\gamma
%continuous. (Cf. Theorem 26.6.viii of Kallenberg.)
%By Kallenberg's version of Knight's theorem for stable integrals (see
%process, say $\tilde A$, which has the same law as $A$.
%% No es exactamente el teorema de Kallenberg... lo que debo hacer es
%encontrar una caracterizaci\'on de procesos estables. Tal vez si
%calculo la variaci\'on cuadr\'atica previsible y la exponencial estoc
Then, we can use Kallenberg's results on time-changes of stable
stochastic integrals found in~\cite{MR1158024} (cf. equation (1.4)) to
infer first the existence of a L\'evy process $\tilde A$ with the same
law as $A$ such that for $t\geq0$:
\[
A_{\mathrm{e}^{\alpha\gamma t}}-A_{1}\stackrel{d} {=}A_{\mathrm{e}^{\alpha\gamma
t}-1}=A_{\int_0^t ( \mathrm{e}^{ \gamma s}/C) ^{\alpha} \,\mathrm{d}s}=
\frac
{1}{C}\int_0^t \mathrm{e}^{\gamma s} \,\mathrm{d}
\tilde A_s.
\]
Then, by associativity of the stochastic integral, we see that
\[
\int_0^t Ce^{-\gamma s}
\,\mathrm{d}A_{\mathrm{e}^{\gamma s \alpha}}=\tilde A_t.
\]
(A similar argument extends when integrating from $s$ to $t$.) Hence,
$\tilde X$ is a stationary version of the Ornstein--Uhlenbeck process
driven by $A$.

On the other hand, the zero set of $\tilde X$ is then the logarithm of
the zero set of $A$. The latter is a self-similar regenerative set and
therefore a random cutout set. Actually, the zero set of $A$ is empty
if $\alpha\in(0,1]$ (cf. \cite{Ber96}, page 63) while if $\alpha\in
(1,2]$ then it has the law of the (closure) of the image of a $\beta
$-stable subordinator with $\beta=1-1/\alpha$ (see \cite{MR0138134} and~\cite{MR0145596}).
%(Still unsure about $\beta$, a place to start is \cite{MR0145596}. The
%asymptotic distribution of the number of zero-free intervals of a
%stable process by Getoor.)

Let
\[
\Xi=\sum_{i\in\mathcal{I}}\delta_{(t_i,x_i)}
\]
be a Poisson point process with intensity $\mathrm{d}t\otimes\mu$ where $\mu
(\mathrm{d}x)=(1-\beta) x^{-2} \,\mathrm{d}x$. Then the zero set of $A$ has the law of a
random cutout set based on $\Xi$ and so the zero set of $\tilde X$ is
the random cutout set (on $\mathbb{R}$) obtained by removing the intervals
\[
(s_i, s_i+z_i)= \bigl(\log(t_i),
\log(t_i+x_i) \bigr).
\]
Namely, we have
\[
\mathcal{Z}:=\{t\in\mathbb{R}; \tilde{X}_{t}=0\}=\mathbb{R}- \bigcup
_{i\in I}\,]s_{i},s_{i}+z_{i}[.
\]
The intensity of the point process
\[
\tilde\Xi=\sum_{i\in\mathcal{I}} \delta_{(s_i,z_i)}
\]
is the image of the measure $\mathrm{d}t\otimes\mu(\mathrm{d}x)$ by $(t,x)\mapsto
(s,z)=(\log(t), \log(t+x))$. A notable cancellation occurs and it is
found to be equal to $(1-\beta)\mathrm{e}^{z}/(\mathrm{e}^{z}-1)^2 \,\mathrm{d}s \,\mathrm{d}z$.
%$(1-\beta)\mathrm{e}^{\tilde x}/(\mathrm{e}^{\tilde x}-1)^2 \,\mathrm{d}\tilde t d\tilde x$.
Using the identity established between the cutting measures of random
cutouts on $(-\infty,\infty)$ and random cutouts on $(0,\infty)$ (cf.
Theorem 2 of \cite{MR799145}) %Fitzsimmons-Fristedt-Shepp
we see that the zero set of the $\alpha$-stable Ornstein--Uhlenbeck
process (started at zero) is a random cutout set whose cutting measure
has density
\[
(1-\beta)\mathrm{e}^{z}/ \bigl(\mathrm{e}^{z}-1 \bigr)^2\,\mathrm{d}z.
%(\tilde t,\tilde x)\mapsto(1-\beta)\mathrm{e}^{\tilde x}/(\mathrm{e}^{\tilde x}-1)^2 d
\]
\upqed\end{pf}
The cutting measure with density $(1-\beta)\mathrm{e}^{z}/(\mathrm{e}^{z}-1)^2$ was
studied in detail in Example 8 of~\cite{MR799145}.
%$(1-\beta)\mathrm{e}^{\tilde x}/(\mathrm{e}^{\tilde x}-1)^2$
There, the authors show that the associated subordinators have zero
drift and L\'evy measure $\nu$ given by
\[
\nu (x,\infty) =\frac{C}{(\mathrm{e}^x-1)^{1-\beta}}.
\]
Note that the density of $\nu$ can also be written as
\[
\nu(\mathrm{d}x)=\frac{C \mathrm{e}^x}{(\mathrm{e}^x-1)^{1-\beta}} \,\mathrm{d}x=\frac{C' \mathrm{e}^{({\beta
}/{2})x}}{(\sinh(x/2))^{2-\beta}} \,\mathrm{d}x.
\]
Hence, in the special case of the Ornstein--Uhlenbeck process associated
to Brownian motion, we recover the results of \cite{MR1478737} (which
go back to \cite{MR1166408}). However, we also deduce that in the
general stable case, the zero set is the image of the Lamperti stable
subordinators introduced in \cite{MR2341013} and studied in general in
\cite{MR2792485}.
%
%co15 #&#
\begin{cor} If $\alpha\in(1,2]$, the random set $\mathcal
{Z}=\overline
{\{t\geq0, X_{t}=0\}}$ is infinitely divisible. Moreover, $\mathcal
{Z}$ is almost surely not bounded ($0$ is recurrent) and we have for
all $t>0$
\[
\operatorname{dim}_{H} \bigl(\mathcal{Z}\cap[0,t] \bigr)=1/\alpha.
\]
\end{cor}
\begin{pf} We only have to give a proof for the Hausdorff dimension.
The Laplace exponent of the Lamperti stable subordinator involved in
Theorem \ref{stableOU} is $\kappa(\gamma)=\frac{\Gamma(1-\beta
+\gamma
)}{\Gamma(1-\beta)\Gamma(\gamma)}$ (see equation (27) in~\cite
{MR799145}). Therefore $\operatorname{dim}_{H}(\mathcal{Z}\cap
[0,t])=\underline
{\operatorname{Ind}}(\kappa)=1-\beta=1/\alpha$.
\end{pf}
\section*{Acknowledgments} C.~Foucart would like to express his gratitude to
Jean Bertoin and thanks Xan Duhalde for pointing the equalities in
Lemma 5.6 of Duquesne and Le Gall's paper \cite{DuqLG}. The authors
would like to thank Arno Siri-J\'egousse for making them aware of each
other's work and thank the referee for his/her careful reading.
%% \MRhref is called by the amsart/book/proc definition of \MR.
% \href{http://www.ams.org/mathscinet-getitem?mr=#1}{#2}
%}
%
% imsref loaded by akundreckaite, 2013-08-13 16:39:41
%

%


%suskaldyti doi

\printhistory

\end{document}